\documentclass[11pt,leqno,twoside]{article}

\usepackage{amsfonts,amsmath,amsthm,amssymb}
\usepackage{enumerate}
\usepackage{graphics,graphicx,subfigure}

\usepackage{color}
\usepackage{todonotes}
\usepackage{cancel}
\usepackage{url}
\usepackage{hyperref}
\usepackage{makeidx}
\usepackage{showidx}
\usepackage{multicol}        
\usepackage{xspace}
\usepackage{stmaryrd}        
\usepackage{pifont}          
\usepackage{fancybox}        
\usepackage{bm}

 \usepackage{fancyhdr}
\theoremstyle{plain}
 \theoremstyle{definition}
 \newtheorem{lem}{Lemma}
 \newtheorem{defn}[lem]{Definition}
 \newtheorem{thm}[lem]{Theorem}
 \newtheorem{prop}[lem]{Proposition}
 \newtheorem{cor}[lem]{Corollary}
 \newtheorem{notn}[lem]{Notations}
 \newtheorem{pb}[lem]{Problem}
 \newtheorem{form}[lem]{Formulae}
 
 \newtheorem*{rk}{Remark}
 \newtheorem*{com}{Comment}
 \newtheorem*{ex}{Example}
 \theoremstyle{remark}

 \newcommand{\blem}{\begin{lem}}
 \newcommand{\elem}{\end{lem}}
 \newcommand{\bdefn}{\begin{defn}}
 \newcommand{\edefn}{\end{defn}}
 \newcommand{\bthm}{\begin{thm} }
 \newcommand{\ethm}{\end{thm}}
 \newcommand{\bprop}{\begin{prop}}
 \newcommand{\eprop}{\end{prop}}
 \newcommand{\bcor}{\begin{cor}}
 \newcommand{\ecor}{\end{cor}}
 \newcommand{\bnotn}{\begin{notn}}
 \newcommand{\enotn}{\end{notn}}
 \newcommand{\bpb}{\begin{pb}}
 \newcommand{\epb}{\end{pb}}
 \newcommand{\bform}{\begin{form}}
 \newcommand{\eform}{\end{form}}
 \newcommand{\brk}{\begin{rk}}
 \newcommand{\erk}{\end{rk}}
 \newcommand{\bcom}{\begin{com}}
 \newcommand{\ecom}{\end{com}}
 \newcommand{\bex}{\begin{ex}}
 \newcommand{\eex}{\end{ex}}
 \newcommand{\bpf}{\begin{proof}}
 \newcommand{\epf}{\end{proof}}

\newcommand{\vb}{{\bf b}}

\newcommand{\ve}{{\bf e}}
\newcommand{\vf}{{\bf f}}
\newcommand{\vg}{{\bf g}}
\newcommand{\vh}{{\bf h}}

\newcommand{\vu}{{\bf u}}
\newcommand{\vv}{{\bf v}}

\newcommand{\cA}{\mathcal{A}}

\newcommand{\cC}{\mathcal{C}}

\newcommand{\cE}{\mathcal{E}}

\newcommand{\cK}{\mathcal{K}}

\newcommand{\cS}{\mathcal{S}}

\newcommand{\cX}{\mathcal{X}}

\newcommand{\bN}{\mathbb{N}}

\newcommand{\bR}{\mathbb{R}}

\newcommand{\be}{\begin{equation}}
\newcommand{\ee}{\end{equation}}
\newcommand{\bal}{\begin{align}}
\newcommand{\eal}{\end{align}}
\newcommand{\ba}{\begin{align*}}
\newcommand{\ea}{\end{align*}}
\newcommand{\bmx}{\begin{matrix}}
\newcommand{\emx}{\end{matrix}}
\newcommand{\bbmx}{\begin{bmatrix}}
\newcommand{\ebmx}{\end{bmatrix}}
\newcommand{\bpmx}{\begin{pmatrix}}
\newcommand{\epmx}{\end{pmatrix}}
\newcommand{\bvmx}{\begin{vmatrix}}
\newcommand{\evmx}{\end{vmatrix}}

\newcommand{\ol}{\overline}
\newcommand{\wh}{\widehat}
\newcommand{\wt}{\widetilde}
\newcommand{\f}{\frac}

\newcommand{\inc}{\subseteq}

\newcommand{\tn}[1]{{|\!|\!|}{#1}{|\!|\!|}}
\newcommand{\tnBig}[1]{{\Big|\!\Big|\!\Big|}{#1}{\Big|\!\Big|\!\Big|}}

\newcommand{\Id}{\mathrm{Id}}

\newcommand{\minimize}[1]{\underset{#1}{\rm minimize}\,}

\newcommand{\la}{\lambda}
\newcommand{\La}{\Lambda}
\newcommand{\eps}{\varepsilon}
  
\title{\vspace{-20mm}Optimal Prediction of Multivalued Functions from Point Samples\medskip\hrule height 1.2pt \vspace{-6mm}}
\author{Simon Foucart\footnote{S. F. is partially supported by  a grant from the Office of Naval Research (N00014-20-1-2787).
Part of this work was carried out during a stay at the Isaac Newton Institute, supported by a grant from the Heilbronn Institute.}  --- Texas A\&M University}
\date{\vspace{-6mm}\rule{100mm}{0.8pt}}

\newcommand\shorttitle{Optimal Prediction of Multivalued Functions from Point Samples}
\newcommand\authors{S.~Foucart}

\begin{document}
\maketitle

\vspace{-15mm}
\begin{abstract}
Predicting the value of a function $f$ at a new point given its values at old points is an ubiquitous scientific endeavor,
somewhat less developed when $f$ produces multiple values that depend on one another,
e.g. when it outputs likelihoods or concentrations.
Considering the points as fixed (not random) entities and focusing on the worst-case,
this article uncovers a prediction procedure that is optimal relatively to some model-set information about $f$.
When the model sets are convex, this procedure turns out to be an affine map constructed by solving a convex optimization program.
The theoretical result is specified in the two practical frameworks of (reproducing kernel) Hilbert spaces and of spaces of continuous functions. 
\end{abstract}



\vspace{-5mm}
\begin{center}
\rule{100mm}{0.8pt}
\end{center}


\section{Introduction}
\vspace{-3mm}

In this article, one considers the general problem of predicting the value $f(x^{(0)})$ of an unknown function $f$ at a `new' point $x^{(0)}$,
given $M$ pieces of data 
$$
y_m = f(x^{(m)}),
\qquad
i \in [1:M],
$$
consisting of evaluations of $f$ at `old' points $x^{(1)},\ldots,x^{(M)}$.
This undertaking is a core problem of Statistical Learning Theory
(see e.g. \cite{ShaBen} for an initiation).
There, the $x^{(m)}$'s are considered independent realizations of a random variable
and the performance of a prediction process is assessed by a notion of risk, 
which reflects a focus on the average case.
The risk decomposes as the sum of an estimation error, often carefully studied, and an approximation error, typically postulated to be small---thus implicitly assuming that the unknown function $f$ is well approximated by elements of a chosen hypothesis class.
Another view, favoring a focus on the  worst case,
is provided by the theory of Optimal Recovery
(see e.g. \cite{MicRiv} for an initiation).
There, the $x^{(m)}$'s are considered fixed deterministic entities and one has some {\sl a priori} knowledge about $f$ in the form $f \in \cK$---for instance, that $f$ is approximated by elements of a hypothesis class up to some given accuracy.
The prediction process, 
encapsulated by a map $\Delta$ taking the $y_m$'s as input and returning an estimation of $f(x^{(0)})$ as an output,
has its performance quantified by the (global\footnote{To distinguished from the local worst-case error, also known as Chebyshev radius.}) worst-case error
\be
\label{GWCE0}
{\rm gwce}(\Delta)  := \sup_{f \in \cK} \big\|f(x^{(0)}) - \Delta (f(x^{(1)}),\ldots, f(x^{(M)})) \big\|.
\ee
A prime objective of Optimal Recovery
is to find a map $\Delta$ that minimizes ${\rm gwce}(\Delta)$, 
or nearly does. 
A standard result in the field asserts that such an optimal $\Delta$ can be chosen as a linear map
as soon as $\cK$ is convex and symmetric about the origin and the evaluation map $f \mapsto f(x^{(0)})$ is a linear functional,
i.e., the function $f$ outputs a single value.

However, in many Data Science settings,
the functions of interest instead output several values.
As an important example,
the functions $\vf$ associated with a neural network are not only multivariate (the inputs are in $\bR^D$)
but also multivalued (the outputs are in $\bR^N$).
When neural networks are used for classification,
the $n$th component $f_n$ of $\vf$ often represents the likelihood of the $n$th class,
so $f_1,\ldots,f_N$ are furthermore nonnegative and sum up to one.
This situation is common in practical problems,
e.g. $\vf(x) = [f_1(x);\ldots;f_N(x)]$ can represent chemical concentrations evolving with time~$x$.
The illustrative theorems stated below,
dealing with such as case, 
instantiate the more general results proved later,
namely Corollaries \ref{CorResHil} and \ref{CorCNonneg}.

The first theorem takes place in the setting of a reproducing kernel Hilbert space $H$ with kernel~$K$, say.
The dependence relation $f_1 + \cdots + f_N = 1$ is imposed,
but not the nonnegativity of each component $f_n \in H$.
These components are assumed to be well approximated by hypothesis classes which are linear subspaces of $H$ with finite dimensions.
Note that the dependence relation presumes that the constant functions belong to $H$,
which excludes the Gaussian kernel, see~\cite{Min}.
Note also that the norm on $\bR^N$ used in \eqref{GWCE0} to quantify the worst-case error is chosen to be the $\ell_\infty$-norm. 

\bthm
\label{ThmIntroHilb}
In a reproducing kernel Hilbert space $H$ containing the constant functions,
consider multivalued functions $\vf \in H^N$ satisfying 
${\rm dist}_H(f_1,V_1) \le \eps_1,\ldots, {\rm dist}_H(f_N,V_N) \le \eps_N$ for some finite-dimensional linear subspaces $V_1,\ldots,V_N$ of $H$ and some parameters $\eps_1,\ldots,\eps_N \ge 0$, as well as $f_1+\cdots+f_N = 1$. 
The prediction of $\vf(x^{(0)}) \in \bR^N$
from the point values $y_1 = \vf(x^{(1)})$ $,\ldots,$ $y_M = \vf(x^{(M)}) \in \bR^N$
is accomplished with minimal $\ell_\infty$-worst-case error
by an affine recovery map $\Delta^{\rm aff} : y \in \bR^{M \times N} \mapsto \sum_{m=1}^M \sum_{n=1}^N y_{m,n} \wh{c}^{(m,n)}  + \wh{d} \in \bR^N$.
Its construction involves,
for each $j \in [1:N]$,
the solutions $\wh{c}_j = \left[ \wh{c}_j^{(m,n)}
\right]_{\substack{m \in [1:M]\\ n \in [1:N] }} \in \bR^{M \times N}$
and $\wh{d}_j \in \bR$ of the convex optimization program \vspace{-2mm}
\begin{align*}
\minimize{\substack{ u^+,u^- \in H \\ c \in \bR^{M \times N}, \, d,e \in \bR }} 
 \;   e 
 & & \mbox{s.to} & \quad 
\pm  \langle u^\pm, 1 \rangle
+ \sum_{n=1}^N \eps_n \,
\Big\| \delta_{j,n} K(\cdot,x^{(0)}) - \sum_{m=1}^M c_{m,n} K(\cdot,x^{(m)}) -  u^\pm  \Big\|_H
\le e \pm d \\
& & \mbox{and} & \quad \;
 \delta_{j,n} v(x^{(0)}) - \sum_{m=1}^M c_{m,n} v(x^{(m)}) - \langle u^\pm, v \rangle = 0
\quad \mbox{for all }v \in V_n, \; n \in [1:N].
\end{align*}
\ethm

This result provides a genuinely practical way of constructing an optimal recovery map,
since the optimization program---a second-order come program---can be solved computationally, albeit not in closed form.
This is true even in the (predominant) case where $H$ is infinite dimensional:
in the spirit of a representer theorem,
the minimizers can be searched for in a finite-dimensional subspace of $H$,
see Subsection \ref{SubSecHil} for details.

The second theorem takes place in a space $C(\cX)$ of continuous functions on a compact set $\cX$.
Still assuming that the components are well approximated by finite-dimensional subspaces and sum up to one,
one also incorporates here the fact that they are nonnegative.
In the statement below, the set $B(\cX)$, resp. $B_+(\cX)$, represents the set of signed, resp. nonnegative, Borel measures on $\cX$.

\bthm
\label{ThmIntroCont}
In a space $C(\cX)$ of continuous functions on a compact set $\cX$,
consider multivalued functions $\vf \in C(\cX)^N$ satisfying 
${\rm dist}_{C(\cX)}(f_1,V_1) \le \eps_1,\ldots, {\rm dist}_{C(\cX)}(f_N,V_N) \le \eps_N$ for some finite-dimensional linear subspaces $V_1,\ldots,V_N$ of $C(\cX)$ and some parameters $\eps_1,\ldots,\eps_N \ge 0$, as well as $f_1+\cdots+f_N = 1$ and $f_1 \ge 0, \ldots, f_N \ge 0$. 
The prediction of $\vf(x^{(0)}) \in \bR^N$
from the point values $y_1 = \vf(x^{(1)}), \ldots, y_M = \vf(x^{(M)}) \in \bR^N$
is accomplished with minimal $\ell_\infty$-worst-case error
by an affine recovery map $\Delta^{\rm aff} : y \in \bR^{M \times N} \mapsto \sum_{m=1}^M \sum_{n=1}^N y_{m,n} \wh{c}^{(m,n)}  + \wh{d} \in \bR^N$.
Its construction involves,
for each $j \in [1:N]$,
the solutions $\wh{c}_j = \left[ \wh{c}_j^{(m,n)}
\right]_{\substack{m \in [1:M]\\ n \in [1:N] }} \in \bR^{M \times N}$
and $\wh{d}_j \in \bR$ of the convex optimization program\vspace{-2mm}
\begin{align*}
\minimize{\substack{ \mu^+, \mu^- \in B(\cX)\\ \bm{\nu}^+, \bm{\nu}^- \in B(\cX)^N \\ c \in \bR^{M \times N}, \, d,e \in \bR }} 
\; \,  e 
& & \mbox{s.to} & \;
\pm \int_{\cX} d\mu^\pm
+ \sum_{n=1}^N \eps_n \int_\cX d|\nu^\pm_n|
\le e \pm d,\\
& & \mbox{and} & \;
\int_\cX v \, d{\nu^\pm_n} = 0 \quad
\mbox{for all } v \in V_n, \; n \in [1:N],\\
& & \mbox{and} & \; \; 
\nu^\pm_n \mp \Big( \delta_{j,n} \delta_{x^{(0)}} - \sum_{m=1}^M c_{m,n} \delta_{x^{(m)}} - \mu^\pm \Big)
\in B_+(\cX),
\quad n \in [1:N].
\end{align*}
\ethm

The genuine practicality of the construction is now a more subtle question,
because the optimization variables include measures, 
which are infinite dimensional objects.
But optimizing over measures can be performed by solving a hierarchy of semidefinite programs,
in the sense that solutions to the infinite-dimensional program are limits of solutions to finite-dimensional semidefinite programs of increasing sizes 
and, in favorable situations,
the convergence occurs in a finite number of steps.
The readers are referred to \cite{JBL} and the references therein for more details on this so-called Moment-SOS hierarchy.
The current situation has an extra advantage,
namely that one is not really interested in the minimizing measures but rather in the minimizing coefficients $c \in \bR^{M \times N}$ and $d \in \bR$.
This point is further discussed in Subsection \ref{SubSecCX}.

Here is an outline of the article's content leading to Theorems \ref{ThmIntroHilb} and \ref{ThmIntroCont} above.
In Section \ref{SecPrel},
after a succinct rundown on Optimal Recovery,
it is proved that affine maps provide optimal procedures for the prediction, and more generally for the estimation of linear functionals, of multiobjects.
This statement covers the case of independent components relative to any $\ell_p$-worst-case error 
and the case of dependent components relative to the $\ell_\infty$-worst-case error.
Section \ref{SecMain} concentrates on the latter case and transforms the existence result into a constructive, yet rather abstract, result 
(Theorem \ref{ThmAbs}) involving the support functions of convex model sets.
The main ingredient is a functional analytic result stemming from the dominated extension theorem and exploited through two of its consequences.
The constructive result is also particularized for two specific model sets.
Next, in Section \ref{SecInst},
the results are applied in the more practical settings of Hilbert spaces and of spaces of continuous functions,
for which the predecessors of Theorems \ref{ThmIntroHilb} and \ref{ThmIntroCont} are derived.
Computability is also discussed in this section.
Section \ref{SecConc} gives an outlook on further aspects to be explored around the recovery of multivalued functions.
The article finishes with an appendix,
of interest in its own right,
where some Optimal Recovery facts are combined and revisited.
\vspace{-2mm}

\paragraph{Notation.}
Throughout this article, 
plain letters, such as $f,g,h$, are used to denote single-valued functions or more generally elements of a vector space $F$.
Boldface letters, such as $\vf, \vg, \vh$,
are used to denote multivalued functions,
or more generally elements in $F^N$,
e.g. $\vf = [f_1; \ldots; f_N]$ with $f_1,\ldots,f_N \in F$.
Indices are written in lowercase,
while their maximal value is written in uppercase,
so e.g. $n \in [1:N]$ means that the index $n$ ranges from $1$ to $N$.
Linear functionals are designated by Greek letters,
such as $\la,\mu,\nu \in F^*$,
where $F^*$ stands for the dual space of $F$.
Elements of $(F^*)^N$ are again written in boldface,
as in $\bm{\mu} = [\mu_1;\ldots;\mu_N]$.
Particular linear functionals are the evaluations at points $x$,
which are symbolized by $\delta_x$,
so that $\delta_x(f) = f(x)$ when $f$ is a function.
There should not be any confusion with the Kronecker symbol $\delta_{i,j}$.
The Greek letter $\rho$ is used for the Minkowski functional $\rho_\cS$ associated with a set $\cS \inc F$.
As for the support function of $\cS$,
is is written as $\tn{\eta}_\cS$ when evaluated at some $\eta \in F^*$.
The dual cone of $\cS$, which is a subset of $F^*$,
is denoted by $\cS^{\rm dual}$.   
Specific spaces $F$ coming into play are Hilbert spaces,
designated by $H$,
and spaces of continuous functions on a compact set $\cX$,
designated by $C(\cX)$.
The set of signed, resp. nonnegative, Borel measures on $\cX$ is denoted as $B(\cX)$, resp. $B_+(\cX)$.

\section{Preliminary Considerations}
\label{SecPrel} 
\vspace{-3mm}

In this section,
after recalling some classical results from Optimal Recovery,
one quickly unravels a solution to the problem of optimally predicting multivalued functions that present no dependence between their components.
Then one establishes the optimality of affine maps 
when such a dependence is integrated as prior knowledge in the set $\cK$, required to be convex,
and when the output space is endowed with the $\ell_\infty$-norm.

\subsection{Brief overview of Optimal Recovery}
\label{SSecOverviewOR}
\vspace{-2mm}

As a start,
it is appropriate to recall some ingredients that are already available in the theory of Optimal Recovery.
In its abstract form,
one does not restrict $f$ to be a function---it generally denotes an element from a vector space $F$,
so multivalued functions are covered.
The {\sl a priori} information reads $f \in \cK$ for some subset $\cK$ of $F$.
This so-called model set encapsulates e.g. scientific knowledge about~$f$.
As for the {\sl a posteriori} information,
it is provided by data $y_1 = \la_1(f),\ldots, y_M  = \la_M(f)$
obtained from linear maps\footnote{The $\la_m$'s are usually linear functionals, so that $Y_m = \bR$, and one can always reduce to this case, but if $f$ is an $N$-valued function and $\la_m$ is a point evaluation, then $Y_m = \bR^N$.} $\la_m: F \to Y_m$ applied to $f$---they do not need to come from point evaluations.
One condenses this information as $y = \La f$,
where the linear map $\La: F \to Y_1 \times \cdots \times Y_M$ is called observation map.
The goal is then to estimate $\Gamma(f)$,
where the quantity of interest $\Gamma: F \to Z$ is a linear map.
In the prediction problem, one takes $\Gamma$ to be
the point evaluation $\delta_{x^{(0)}}: f \in F \mapsto f(x^{(0)}) \in \bR^N$,
but again no restriction on $\Gamma$ other than linearity is imposed.
Thus, the choice $\Gamma = \Id_F$ corresponds to the full recovery problem.
The estimation process boils down to a map 
(cognizant of $\cK$ and $\La$)
$\Delta: Y_1 \times \cdots \times Y_M \to Z$,
whose performance is assessed,
as in \eqref{GWCE0}, by the worst-case error
$$
{\rm gwce}(\Delta)  := \sup_{f \in \cK} \big\|\Gamma(f) - \Delta (\La f) \big\|_Z.
$$
The practical construction of recovery maps $\Delta^{\rm opt}$ that are optimal 
(meaning ${\rm gwce}(\Delta^{\rm opt}) \le {\rm gwce}(\Delta)$ for any~$\Delta$)
or near-optimal (meaning ${\rm gwce}(\Delta^{\rm opt}) \le C \,{\rm gwce}(\Delta)$ for any~$\Delta$)
is the ultimate goal of Optimal Recovery.
One aims at such constructions in as many relevant situations as possible, always with the hope that $\Delta^{\rm opt}$ shall be a `simple' map, 
the epitome of which is a linear map.
Along these lines, the following facts are well established: \vspace{-5mm}
\begin{enumerate}

\item[(i)] If the quantity of interest $\Gamma$ is a linear functional and if the model set $\cK$ is convex and symmetric about the origin,
then there is a linear recovery map which is optimal;
furthermore, if $\cK$ is merely convex, then the conclusion holds with linear replaced by affine
(this fact is due to \cite{Suk} and is reproved in the appendix with a different argument).

\item[(ii)] If $F$ is a Hilbert space and if the model set $\cK$ is a centered hyperellipsoid, then there is a linear recovery map which is optimal 
(this fact is also reproved unconventionally in the appendix);
furthermore, 
if $\cK$ is the intersection of two (not more) centered hyperellipsoids, then the conclusion still holds
(this fact was recently established in \cite{2Hyp}).

\item[(iii)] Independently of $F$ and $\Gamma$,
if the model set $\cK$ is convex and symmetric about the origin,
then there is always a linear recovery map which is near-optimal with factor $C = 1 + \sqrt{M}$
(this fact is from \cite{CreWoj});
furthermore,  if $\cK$ is merely convex, then the conclusion holds for affine maps with a near-optimality factor $C = 2 \min\{ \sqrt{N},  1 + \sqrt{M}\}$
(this fact is proved in the appendix).\vspace{-3mm}
\end{enumerate}

For the prediction of multivalued functions,
none of these results are very useful:
(i) does not apply because, as already pointed out, the point evaluation $\delta_{x^{(0)}}$ maps into $\bR^N$ and hence is not a linear functional,
(ii) is too restrictive on the choice of model set,
and (iii) supplies near-optimality factors which are too large.

\subsection{Simple solution in case of independent components}
\vspace{-2mm}

In this subsection,
the object is not yet pinned down as a multivariate function,
as it is only required to be of the form $\vf = [f_1;\ldots;f_N]$ where each component $f_n$ belongs to a common vector space $F$.
For each $n$,
there is an {\sl a priori} information of the type $f_n \in \cK_n$ for some subset $\cK_n$ of $F$.
Thus, the overall model set for the whole object is
\be
\label{IndepMS}
\bm{\cK}^{\rm ind}: = \big\{
\vf = [f_1;\ldots;f_N] \in F^N: f_1 \in \cK_1, \ldots, f_N \in \cK_N
\big\} \inc F^N.
\ee 
No dependence relation between the $f_n$'s,
such as $f_1 + \cdots + f_N = 1$,
is assumed at this point---this
explains the chosen terminology of independent components.
As soon to be revealed, this situation is not too interesting,
since optimal recovery maps can be constructed componentwise.
In the following statement,
the linear functionals $\la_1,\ldots,\la_M \in F^*$
generalize the point evaluations $\delta_{x^{(1)}},\ldots,\delta_{x^{(M)}}$ and give rise to the observation map
\be
\label{IndObsMap}
\La: \vf \in F^N \mapsto [\la(f_1),\ldots,\la(f_N)] \in \bR^{M \times N},
\qquad
\la(f_n) := [\la_1(f_n);\ldots;\la_M(f_n)] \in \bR^M,
\ee
while the quantity of interest 
\be
\label{IndQoI}
\Gamma : \vf \in F^N \mapsto [\gamma_1(f_1);\ldots;\gamma_N(f_N)] \in \bR^N 
\ee
should be thought of with linear functionals $\gamma_1,\ldots,\gamma_N$ all being  equal to $\delta_{x^{(0)}}$. 
With only minor modifications, the proof below would show that the optimal recovery map $\Delta^{\rm aff}$ can be taken linear if the $\cK_n$'s are not only convex but also symmetric about the origin.

\bprop
\label{PropIndep}
Given the model set \eqref{IndepMS} based on convex sets $\cK_1,\ldots,\cK_N$,
the observation map~\eqref{IndObsMap},
and the quantity of interest~\eqref{IndQoI},
the worst-case error of a recovery map $\Delta: \bR^{M \times N}~\to~\bR^N$,
 measured in $\ell_p$ for $p \in [1,\infty]$ as
$$
{\rm gwce}_p(\Delta)  := \sup_{\vf \in \bm{\cK}^{\rm ind}} \big\|\Gamma(\vf) - \Delta (\La \vf) \big\|_{\ell_p^N},
$$
is minimized for the affine map $\Delta^{\rm aff}: \bR^{M \times N}~\to~\bR^N$ given by
\be
\label{ORM-Ind}
\Delta^{\rm aff}(y) 
= [\Delta^{\rm aff}_1(y_1); \ldots; \Delta^{\rm aff}_N(y_N)] \in \bR^N,
\qquad 
y = [y_1, \ldots, y_N] \in  \bR^{M \times N},
\ee
each $\Delta_n^{\rm aff}: \bR^M \to \bR$ being an affine optimal recovery map for the estimation of $\gamma_n$ relative to~$\cK_n$.   
\eprop

\bpf
A well-known lower bound for the worst-case error of any recovery map $\Delta$ involves the symmetrized set $\bm{\cK}^{\rm ind}_{\rm sym} := (\bm{\cK}^{\rm ind} - \bm{\cK}^{\rm ind})/2$.
It is obtained by considering any $\vh \in \ker(\La) \cap \bm{\cK}^{\rm ind}_{\rm sym}$,
writing $\vh = (\vf^+ - \vf^-)/2$ with $\vf^+,\vf^- \in \bm{\cK}^{\rm ind}$ and $\La(\vf^+) = \La(\vf^-) =: y$, 
and deriving
$$
{\rm gwce}_p(\Delta) \ge 
\f{1}{2} \big\|\Gamma(\vf^+) - \Delta (y) \big\|_{\ell_p^N} 
+ \f{1}{2} \big\|\Gamma(\vf^-) - \Delta (y) \big\|_{\ell_p^N}
\ge \f{1}{2} \big\|\Gamma(\vf^+) - \Gamma(\vf^-) \big\|_{\ell_p^N}
= \big\| \Gamma(\vh) \big\|_{\ell_p^N}.
$$
Taking the supremum over $\vh$ results in the said lower bound, which reads 
$$
{\rm gwce}_p(\Delta) \ge \sup_{\substack{\vh \in \ker(\La) \\ \vh \in \bm{\cK}^{\rm ind}_{\rm sym} }} \|\Gamma(\vh)\|_{\ell_p^N}.
$$
Noticing
that $\vh \in \ker(\La)$ if and only if $h_1 \in \ker(\la), \ldots, h_N \in \ker(\la)$,
that
$\vh \in \bm{\cK}^{\rm ind}_{\rm sym}$ if and only if $h_1 \in (\cK_1 - \cK_1)/2, \ldots, h_N \in (\cK_N - \cK_N)/2$,
and that $\|\Gamma(\vh)\|_{\ell_p^N}^p = \sum_{n=1}^N |\gamma_n(h_n)|^p$, one arrives at
$$
{\rm gwce}_p(\Delta)^p \ge
\sup_{\substack{ h_n \in \ker(\la), \; {\rm all} \; n  \\ h_n \in (\cK_n - \cK_n)/2,  \; {\rm  all} \; n}}
\sum_{n=1}^N |\gamma_n(h_n)|^p
= \sum_{n=1}^N \sup_{\substack{ h_n \in \ker(\la) \\ h_n \in (\cK_n - \cK_n)/2}} |\gamma_n(h_n)|^p.
$$
According to (i-b) of Proposition \ref{PropApp} and its proof 
(the argument behind the fact that convex model sets ensures optimality of affine maps, which is found in the appendix ),
each summand---power $p$ omitted---coincides with the worst-case error of an affine optimal recovery map $\Delta_n^{\rm aff}$
for the estimation of $\gamma_n$ relative to $\cK_n$.
It follows that 
\begin{align*}
{\rm gwce}_p(\Delta)^p 
& \ge \sum_{n=1}^N \sup_{f_n \in \cK_n} |\gamma_n(f_n) - \Delta^{\rm aff}_n (\la(f_n))|^p 
= \sup_{f_n \in \cK_n, \; {\rm  all} \; n} \sum_{n=1}^N  |\gamma_n(f_n) - \Delta^{\rm aff}_n (\la (f_n) ) |^p\\
& = \sup_{\vf \in \bm{\cK}^{\rm ind}} \|\Gamma(\vf) - \Delta^{\rm aff}(\La \vf)\|_{\ell_p^N}^p
= {\rm gwce}_p(\Delta^{\rm aff})^p.
\end{align*}
This inequality establishes the announced optimality of $\Delta^{\rm aff}$ defined in \eqref{ORM-Ind}. 
Of note, the input of each component $\Delta^{\rm aff}_n$ consists of observations about $f_n$ alone,
not about $f_1,\ldots,f_{n-1},f_{n+1},\ldots,f_N$---as one would have anticipated from the independence of $f_1,\ldots,f_N$.
\epf


\subsection{Existence result in case of dependent components} 
\vspace{-2mm}

The case of independent components having been brushed aside,
this subsection is concerned with a model set that includes dependence between the components of a multiobject $\vf \in F^N$---again it need not be pinned down as a multivalued function yet. 
Generalizing the relation $f_1+\cdots+f_N = 1$,
one introduces several dependence relations of the type
$$
a_{\ell,1} f_1 + \cdots + a_{\ell,N} f_N = b_\ell,
\qquad
\ell \in [1:L],
$$ 
which are summarized as $A \vf = \vb$ for fixed $A \in \bR^{L \times N}$ and $\vb \in F^L$.
As such, the contemplated model set takes the form
\be
\label{DepMS}
\bm{\cK}^{\rm dep} := \big\{ \vf = [f_1;\ldots;f_N] \in F^N : f_1 \in \cK_1, \ldots, f_N \in \cK_N, A \vf = \vb \big\} \inc F^N,
\ee
or in short $\bm{\cK}^{\rm dep} = \bm{\cK}^{\rm ind} \cap A^{-1}(\{\vb\})$.
Note that this model set is convex as soon as the sets $\cK_1,\ldots,\cK_N$ are themselves convex.
According to (i) of Proposition \ref{PropApp} in the appendix,
with the same setting as Proposition \ref{PropIndep},
one can already guarantee that affine maps are near optimal with a factor $2 \sqrt{N}$ (and even better when $p>2$).
For $p=\infty$, genuine optimality can actually be achieved---although not constructively at this point---as established below.

\bprop
\label{PropDep}
Given the model set \eqref{DepMS} based on convex sets $\cK_1,\ldots,\cK_N$,
the observation map~\eqref{IndObsMap},
and the quantity of interest~\eqref{IndQoI},
the worst-case error of a recovery map $\Delta: \bR^{ M \times N}~\to~\bR^N$, measured in $\ell_\infty$, 
is minimized for an affine recovery map.
\eprop

\bpf
The worst-case error of any recovery map $\Delta = [\Delta_1;\ldots;\Delta_N]$ is
\begin{align*}
{\rm gwce}_\infty(\Delta) 
& = \sup_{\vf \in \bm{\cK}^{\rm dep}} \|\Gamma(\vf) - \Delta(\La \vf)\|_{\ell_\infty^N}
= \sup_{\vf \in \bm{\cK}^{\rm dep}} \max_{n \in [1:N]} |\gamma_n(f_n) - \Delta_n(\La \vf)|\\
& = \max_{n \in [1:N]} \sup_{\vf \in \bm{\cK}^{\rm dep}}  |\gamma_n(f_n) - \Delta_n(\La \vf)|
= \max_{n \in [1:N] } {\rm gwce}(\Delta_n).
\end{align*}
If $\Delta_n^{\rm aff}: \bR^{M \times N}~\to~\bR$ denotes an optimal recovery map for the estimation of the linear functional $\vf \in F^N \mapsto \gamma_n(f_n) \in \bR$ relative to the convex set $\bm{\cK}^{\rm dep}$ as a whole,
recalling that $\Delta_n^{\rm aff}$ can be taken as an affine map,
then one has 
$$
{\rm gwce}_\infty(\Delta)  \ge  \max_{n \in [1:N] } {\rm gwce}(\Delta_n^{\rm aff})
= {\rm gwce}_\infty(\Delta^{\rm aff}), 
$$
where $\Delta^{\rm aff}: \bR^{M \times N}~\to~\bR^N$ is simply defined as $\Delta^{\rm aff} := [\Delta^{\rm aff}_1; \ldots; \Delta^{\rm aff}_N]$.
Since the latter is an affine map, 
the result is now justified.
But, contrary to Proposition~\ref{PropIndep},
each $\Delta_n^{\rm aff}$ is defined on all of~$\bR^{M \times N}$,
i.e., it~could take inputs obtained from observing components other than $f_n$.
\epf

\section{Constructive Results in Case of Dependent Components}
\label{SecMain}
\vspace{-3mm}

A constructive version of the above result is derived 
in this section,
specifically in Subsection~\ref{SubSecAbs}.
Before that,
one isolates in Subsection~\ref{SubSecLem}
two useful observations:
one is needed for the proof of the main theorem
and the other for the proof of one of two consequences presented in Subsection~\ref{SubSecTwoPart}.
The results throughout this section are still rather abstract---they will be instantiated to more practical situations in the next section.

\subsection{The key lemmas}
\label{SubSecLem}
\vspace{-2mm}

The upcoming Theorem \ref{ThmAbs} and Theorem \ref{ThmIntersection} rely on Lemma \ref{Lem1} and Lemma \ref{Lem2}, respectively.
These two lemmas, stated right below,
are in fact byproducts of an encompassing theorem,
so their proofs are deferred after the statement and justification of the theorem.

\blem
\label{Lem1}
Given a convex subset $\cC$ of a vector space $F$,
a linear map $\cA$ from $F$ into another vector space $G$,
an element $b \in G$,
and a linear functional $\eta \in F^*$,
one has
$$
\sup_{f \in F} \Big\{ 
\eta(f): f \in \cC \mbox{ and } \cA(f) = b \Big\}
= \min_{\mu \in G^*} 
\Big\{ \mu(b) + \sup_{f \in \cC} \; (\eta - \mu \circ \cA)( f) \Big\}.
$$
\elem

\blem
\label{Lem2}
Given a convex set $\cK$ and a convex cone $\cC$ in a vector space $F$ and given a linear functional $\eta \in F^*$, one has
$$
\sup_{f \in F} \Big\{ \eta(f): f \in \cK \cap \cC \Big\}
= \min_{\nu \in F^*} \Big\{ \sup_{f \in \cK} \nu(f): \; \nu - \eta \in \cC^{\rm dual} \Big\},
$$
where $\cC^{\rm dual} := \{ \mu \in F^*: \mu(f) \ge 0 \mbox{ for all } f \in F\}$ denotes the dual cone of $\cC$.
\elem

The encompassing functional analytic result behind these lemmas  is stated as Theorem \ref{ThmMainDuality} below.
It is likely known---the arguments are certainly standard---but a proof is included for completeness.

\bthm
\label{ThmMainDuality}
Given convex subsets $\cC_1,\ldots,\cC_K$ of a vector space $F$
and a linear functional $\eta \in F^*$,
one has
\be
\label{MainDuality}
\sup_{f \in F} \Big\{ \eta(f) : \; f \in \bigcap\nolimits_{k=1}^K \cC_k \Big\} =
\min_{\mu_1,\ldots,\mu_K \in F^*} \Big\{ 
\sum\nolimits_{k=1}^K \sup_{f_k \in \cC_k} \mu_k(f_k) 
 : 
\sum\nolimits_{k=1}^K \mu_k = \eta
\Big\} .
\ee
\ethm

\bpf
Let ${\rm lhs}$ denote the supremum on the left-hand side of \eqref{MainDuality} and let ${\rm rhs}$ denote the infimum (proved to be a minimum later) on its right-hand side.
The inequality ${\rm rhs} \ge {\rm lhs}$ is easily verified.
Indeed, it suffices to observe that,
for feasible $\mu_1,\ldots,\mu_K \in F^*$ and $f \in F$,
$$
\sum\nolimits_{k=1}^K \sup_{f_k \in \cK} \mu_k(f_k) 
\ge \sum\nolimits_{k=1}^K  \mu_k(f) = \eta(f).
$$
The reverse inequality requires more work.
One starts by recalling that the Minkowski functional (aka gauge) associated with a subset $\cC$ of $F$ is defined by $\rho_\cC(f) := \inf\{ t > 0: f \in t \, \cC \} \in [0,+\infty]$ for any $f \in F$.
Recall also that,
if $\cC$ is convex and $0 \in \cC$,
then $\rho_\cC$ is a sublinear functional, 
meaning that $\rho_\cC(t  f) = t \, \rho_\cC(f)$ and $\rho_\cC(f+g) \le \rho_\cC(f)+\rho_\cC(g)$ for all $t>0$ and all $f,g \in F$.
In particular,
selecting an arbitrary $\ol{f} \in \cC_1 \cap \cdots \cap \cC_K$,
the Minkowski functional associated with the set $\cC:= (\cC_1 - \ol{f}) \cap \cdots \cap (\cC_K - \ol{f})$ is sublinear.
The latter is (easily verified to be) given,
for any $f \in F$, by
$\rho_\cC(f) = \max\{ \rho_{\cC_1 - \ol{f}}(f), \ldots, \rho_{\cC_K - \ol{f}}(f) \}$.
With ${\rm lhs}$ still denoting the supremum on the left-hand side of \eqref{MainDuality},
one notices that $\eta(\ol{f}+h) \le {\rm lhs}$,
or $\eta(h) \le {\rm lhs} - \eta(\ol{f})$,
for all $h \in \cC$.
This implies that $\eta(h) \le [{\rm lhs}-\eta(\ol{f})] \, \rho_\cC(h)$ for all $h \in H$.
Therefore,
the linear functional $\nu$ defined on the subspace $S = \{ [h;\ldots;h], h \in F \}$ of $F^K$
by $\nu([h;\ldots;h]) = \eta(h)$
is dominated on $S$ by the sublinear functional $\rho$ defined on $F^K$ by $\rho([h_1;\ldots;h_K]) = [{\rm lhs}-\eta(\ol{f})] \, \max\{ \rho_{\cC_1 - \ol{f}}(h_1), \ldots, \rho_{\cC_K - \ol{f}}(h_K) \}$.
The dominated extension theorem
now ensures the existence of a linear functional $\mu$ defined on $F^K$ such that $\mu_{|S} = \nu$ and $\mu \le \rho$.  
Such a linear functional takes the form $\mu([h_1;\ldots;h_K]) = \mu_1(h_1) + \cdots + \mu_K(h_K)$ for some linear functionals $\mu_1,\ldots,\mu_K \in F^*$.
The fact that $\mu_{|S} = \nu$ implies that $\mu_1 + \cdots + \mu_K = \eta$.
As for the fact that $\mu \le \rho$,
it implies that 
$\mu_1(h_1) + \cdots + \mu_K(h_K) \le [{\rm lhs}-\eta(\ol{f})] \, \max\{ \rho_{\cC_1 - \ol{f}}(h_1), \ldots, \rho_{\cC_K - \ol{f}}(h_K) \}$ for all $[h_1;\ldots;h_K] \in F^K$.
Thus, given $f_k \in \cC_k$ for $k \in [1:K]$,
since $h_k := f_k - \ol{f} \in \cC_k - \ol{f}$ satisfies $\rho_{\cC_k - \ol{f}}(h_k) \le 1$,
one obtains
$ \mu_1(f_1-\ol{f}) + \cdots + \mu_K(f_K-\ol{f}) \le [{\rm lhs}-\eta(\ol{f})] $.
After simplifying and taking the suprema over $f_1, \ldots, f_K$,
one arrives at 
$$
\sup_{f_1 \in \cC_1} \mu_1(f_1) + \cdots + \sup_{f_K \in \cC_K} \mu_K(f_K)
\le {\rm lhs},
$$
which is the inequality required to complete the proof.
\epf

With Theorem \ref{ThmMainDuality} now justified, 
it remains to deduce Lemmas \ref{Lem1} and \ref{Lem2} as consequences.

\bpf[Proof of Lemma \ref{Lem1}]
Applying Theorem \ref{ThmMainDuality} with $K=2$,
$\cC_1 = \cA^{-1}(\{b\})$,
and $\cC_2 = \cC$ shows that the sought-after supremum equals
$$
\min_{\nu \in F^*}
\bigg\{  \sup_{f_1 \in \cA^{-1}(\{b\})} \nu(f_1) 
+ \sup_{f_2 \in \cC} (\eta - \nu)(f_2)
\bigg\}.
$$
Selecting some $\ol{f} \in \cA^{-1}(\{b\})$,
one observes that
$$
\sup_{f_1 \in \cA^{-1}(\{b\})}  \nu(f_1)
= \sup_{h \in \ker(\cA)} \nu(\ol{f}+h)
= \begin{cases}
\nu(\ol{f}) & \text{if } \ker(\cA) \inc \ker(\nu),\\
+\infty & \text{otherwise}.
\end{cases}
$$
Thus,
the minimum over $\nu \in F^*$ will be achieved by imposing $\ker(\cA) \inc \ker(\nu)$,
which is equivalent to the existence of $\mu \in G^*$ such that $\nu = \mu \circ \cA$,
in which case $\nu(\ol{f}) = \mu(b)$.
The announced result immediately follows.
\epf

\bpf[Proof of Lemma \ref{Lem2}]
Again applying Theorem \ref{ThmMainDuality} with $K=2$, and this time with $\cC_1 = \cK$ and $\cC_2 = \cC$, shows that the sought-after supremum equals
$$
\min_{\nu \in F^*} \bigg\{ \sup_{f_1 \in \cK}\nu(f_1) + \sup_{f_2 \in \cC} (\eta-\nu)(f_2) \bigg\}.
$$
Using the fact that $\cC$ is a cone,
i.e. that if $f \in \cC$,
then $t f \in \cC$ for all $t > 0$,
one observes that 
$$
\sup_{f_2 \in \cC} (\eta-\nu)(f_2) = \begin{cases}
0 & \text{if } -(\eta-\nu) \in \cC^{\rm dual},\\
+\infty & \text{otherwise}.
\end{cases}
$$
Thus,
the minimum over $\nu \in F^*$ will be achieved by imposing $\nu - \eta  \in \cC^{\rm dual}$.
The announced result immediately follows.
\epf

\subsection{Abstract formulation of the main result}
\label{SubSecAbs}
\vspace{-2mm}

It is now time to state the main theorem about optimal recovery of multiobjects whose dependent components belong to convex model sets.
Its formulation involves the support function of a set $\cC \inc F$, as defined by
$$
\tn{\eta}_\cC = \sup_{f \in \cC} \, \eta(f),
\qquad \eta \in F^*. 
$$
Beware that the notation may be misleading:
the support function is not necessarily a norm,
as $\tn{ -\eta }_\cC = \tn{ \eta }_\cC$ might not hold in general.
Nonetheless, it is a sublinear functional---indeed, it is straightforward to verify that it coincides with the Minkowski functional associated with the polar $\cC^\circ := \{ \eta \in F^*: \eta(f) \le 1 \mbox{ for all } f \in \cC \}$ of $\cC$,
which is always convex and contains the origin.
Despite its abstract form, the result is constructive,
in the sense that an affine optimal recovery map is obtained after solving several convex optimization programs.
Note that formally taking $L=0$  (hence discarding $\bm{\mu}^+$ and $\bm{\mu}^-$) also gives a constructive version of Proposition \ref{PropIndep} for the case of independent components.

\bthm
\label{ThmAbs}
In the setting of Proposition \ref{PropDep},
one considers for each $j \in [1:N]$
some solutions $\wh{c}_j = \bbmx \wh{c}^{(m,n)}_j \ebmx_{\substack{ m=1,\ldots,M \\ n=1,\ldots,N }} \in \bR^{M \times N}$ 
and $\wh{d}_j \in \bR$ to the convex optimization program\vspace{-2mm}
\be
\label{MainOptAbs}
\minimize{\substack{ \bm{\mu}^+, \bm{\mu}^- \in (F^*)^L \\ c \in \bR^{M \times N}, \, d,e \in \bR }} 
 \; \,  e \qquad 
\mbox{s.to} \quad
\sum_{\ell = 1}^L 
\mu^\pm_\ell(b_\ell)
+ \sum_{n=1}^N 
\tnBig{ \pm \delta_{j,n} \gamma_j \mp \sum_{m=1}^M c_{m,n} \la_m - \sum_{\ell=1}^L a_{\ell,n} \mu^\pm_\ell 
}_{\cK_n}
\le e \pm d.
\ee
Then an optimal recovery map is given by the affine map
\be
\label{MainAffORMap}
\Delta^{\rm aff} : y \in \bR^{M \times N} \mapsto \sum_{m=1}^M \sum_{n=1}^N y_{m,n} \wh{c}^{(m,n)}  + \wh{d} \in \bR^N.
\ee 
\ethm

\bpf
According to Proposition \ref{PropDep} (and its proof),
one can look for an optimal recovery map 
in the form $\Delta^{\rm aff} = [\Delta^{\rm aff}_1; \ldots; \Delta^{\rm aff}_N]$.
Fixing $j\in [1:N]$, the affine map $\Delta^{\rm aff}_j: \bR^{ M \times N} \to \bR$ can be written as
$$
\Delta^{\rm aff}_j(y) = \sum_{m=1}^M \sum_{n=1}^N y_{m,n} c_{m,n}  + d
$$ 
for some $c \in \bR^{M \times N}$ and $d \in \bR$.
These coefficients are sought to minimize the worst-case error  
$$
{\rm gwce}(\Delta_j^{\rm aff})
= \sup_{\vf \in \bm{\cK}^{\rm dep}} \Big| \gamma_j(f_j ) - \sum_{m=1}^M  \sum_{n=1}^N c_{m,n} \la_m(f_n) - d \Big|
= \sup_{\vf \in \bm{\cK}^{\rm dep}}  \big| \bm{\eta}^{c}(\vf) - d \big|,
$$
with a linear functional $\bm{\eta}^{c}$ implicitly defined on $F^N$ by $\bm{\eta}^{c}(\vf) = \gamma_j(f_j ) - \sum_{m=1}^M  \sum_{n=1}^N c_{m,n} \la_m(f_n)$.
At this point, noticing that
\begin{align}
\nonumber
\sup_{\vf \in \bm{\cK}^{\rm dep}}  \big| \bm{\eta}^{c}(\vf) - d \big|
& = \inf_{e \in \bR} \; e \quad \mbox{s.to } -e \le \bm{\eta}^c(\vf) - d \le e \mbox{ for all } \vf \in \bm{\cK}^{\rm dep}\\
\label{AuxMainRes1}
& = \inf_{e \in \bR} \; e \quad \mbox{s.to }
\begin{cases}
\displaystyle{\sup_{\vf \in \bm{\cK}^{\rm dep}} (+\bm{\eta}^c)(\vf) \le e+d},\\
\displaystyle{\sup_{\vf \in \bm{\cK}^{\rm dep}} (-\bm{\eta}^c)(\vf) \le e-d},
\end{cases}
\end{align}
it remains to invoke Lemma \ref{Lem1} to transform these two suprema into minima over $\bm{\mu}^+,\bm{\mu}^- \in (F^L)^*$.
More precisely,
the first supremum is
$$
\sup_{\vf \in F^N} \Big\{ (+\bm{\eta}^c)(\vf): \vf \in \bm{\cK}^{\rm ind} \mbox{ and } A \vf = \vb \Big\}
= \min_{\bm{\mu}^+ \in (F^L)^*} \Big\{
\bm{\mu}^+(\vb) + \sup_{\vf \in \bm{\cK}^{\rm ind}} (+\bm{\eta}^c - \bm{\mu}^+ \circ A)\vf
\Big\}.
$$
The linear functional $\bm{\mu}^+ \in (F^L)^*$ decomposes through $\mu^+_1,\ldots,\mu^+_L \in F^*$
as $\bm{\mu}(\vb) = \sum_{\ell=1}^L \mu^+_\ell(b_\ell)$,
so $\bm{\mu}^+(A \vf) = \sum_{\ell=1}^L \mu^+_\ell \Big( \sum_{n=1}^N a_{\ell,n} f_n \Big)
= \sum_{n=1}^N  \Big( \sum_{\ell=1}^L a_{\ell,n} \mu^+_\ell \Big)(f_n)$,
and similarly the linear functional $\bm{\eta}^c \in (F^N)^*$ decomposes through $\eta^c_1,\ldots,\eta^c_N \in F^*$ as
$\bm{\eta}^c(\vf) = \sum_{n=1}^N \eta^c_n(f_n)$.
Thus, the first constraint in \eqref{AuxMainRes1}---the one with the $+$ sign---reads:
there exist $\mu^+_1,\ldots,\mu^+_L \in F^*$ such that
$$
\sum_{\ell=1}^L \mu^+_\ell(b_\ell) 
+ \sup_{f_1 \in \cK_1, \ldots, f_N \in \cK_N} 
\sum_{n=1}^N \Big( + \eta^c_n  - \sum_{\ell=1}^L a_{\ell,n} \mu^+_n \Big)(f_n)
\le e + d,
$$ 
in other words, taking the expression of $\eta^c_n$ into account,
$$
\sum_{\ell=1}^L \mu^+_\ell(b_\ell) 
+ \sum_{n=1}^N \sup_{f_n \in \cK_n} \Big( + \delta_{j,n} \gamma_j - \sum_{m=1}^M c_{m,n} \la_m  - \sum_{\ell=1}^L a_{\ell,n} \mu^+_n \Big)(f_n)
\le e + d.
$$ 
This is the constraint appearing in \eqref{MainOptAbs} for the choice $\pm = +$.
Likewise,
the second constraint in \eqref{AuxMainRes1} reduces to the constraint appearing in \eqref{MainOptAbs} for the choice $\pm = -$.
Incorporating the $\mu^\pm_\ell$'s as optimization variables exposes ${\rm gwce}(\Delta_j^{\rm aff})$ as the minimal value of a convex program.
Further minimizing over $c \in \bR^{M \times N}$ and $d \in \bR$ leads to the convex program \eqref{MainOptAbs}. 
Its minimizers provide an optimal recovery map via \eqref{MainAffORMap},
but keep also in mind that its minimal value is equal to the minimal worst-case error over all recovery maps.
\epf

\subsection{Two particular cases of the main result}
\label{SubSecTwoPart}
\vspace{-2mm}

The practical solvability of the convex optimization program \eqref{MainOptAbs} depends on the amenability to computations of the support functions of the model sets $\cK_1,\ldots,\cK_N$.
As such, it is worth looking at two types of particularly relevant sets,
still without specifying the vector space $F$ at this point.

{\bf Approximability sets.}
As alluded to in the introduction, 
model sets based on approximation capabilities are quite pertinent,
since some assumptions about the approximation power of certain hypothesis classes are often made implicitly.
Thus,
it is natural to consider model sets of the type $
\{ f \in F: {\rm dist}_F(f,V) \le \eps\}$
relative to a linear subspace $V$ of $F$ and a parameter $\eps \ge 0 $.
Actually, one considers a more general model set involving an invertible operator $T$ on $F$,
namely
\be
\label{GenApproxSets}
\cK_{V,T} =
\{ f \in F: {\rm dist}_F(Tf,V) \le 1\}.
\ee
Applying Theorem \ref{ThmAbs} then entails to determining the support function of such a set, which is done as follows:
\begin{align}
\nonumber
\tn{\eta}_{\cK_{V,T}}
& = \sup_{f \in F} \{ \eta(f): \mbox{there exists }v \in V \mbox{ such that } \|Tf-v\|_F \le 1 \}\\
\nonumber
&  = \sup_{\substack{f \in F \\ v \in V}} \{ \eta \circ T^{-1} ( Tf) : \| Tf - v\|_F \le \eps \}
 = \sup_{\substack{g \in F \\ v \in V}} \{ \eta\circ T^{-1} (g) + \eta \circ T^{-1} (v)  : \| g\|_F \le 1 \}\\
\label{SupFunAppSet} 
& = \sup_{g \in F } \{ \eta \circ T^{-1} (g)  : \| g\|_F \le 1 \}
+ \sup_{ v \in V} \{  \eta \circ T^{-1} (v)  \}
= \begin{cases}
\|\eta \circ T^{-1}\|_{F^*} & \text{if } (\eta \circ T^{-1} )_{|V} = 0,\\
+\infty & \text{otherwise}.
\end{cases}
\end{align}
This observation leads to the following consequence of Theorem \ref{ThmAbs}
(where the change $\bm{\mu}^- \leftrightarrow -\bm{\mu}^-$ of optimization variable is applied). 
It should be viewed as a generalization to $L \ge 1$ and $T \not= (1/\eps)\Id_F$ of \cite[Theorem 3.1]{DFPW}.

\bthm
\label{ThmAbs4AppSet}
For model sets $\cK_n = \cK_{V_n,T_n}$,
the affine optimal recovery map of Theorem~\ref{ThmAbs} is constructed by solving, for each $j \in [1:N]$,
the constrained convex optimization program \vspace{-2mm}
\begin{align*}
\minimize{\substack{ \bm{\mu}^+,\bm{\mu}^- \in (F^*)^L \\ c \in \bR^{M \times N}, \, d,e \in \bR }} 
\; \, e 
& & \mbox{s.to} & \quad
\pm \sum_{\ell = 1}^L 
\mu^\pm_\ell(b_\ell)
+ \sum_{n=1}^N 
\Big\| \Big(\delta_{j,n} \gamma_j - \sum_{m=1}^M c_{m,n} \la_m - \sum_{\ell=1}^L a_{\ell,n} \mu^\pm_\ell \Big) \circ T_n^{-1} \Big\|_{F^*}
\le e \pm d \\
& & \mbox{and} & \quad
\Big(\delta_{j,n} \gamma_j - \sum_{m=1}^N c_{m,n} \la_m - \sum_{\ell=1}^L a_{\ell,n} \mu^\pm_\ell \Big)(T_n^{-1}(v)) =  0 
\; \, \mbox{for all }v \in V_n, n \in [1:N].
\end{align*}
\ethm


{\bf Model sets intersected with convex cones.}
The above result applies to the recovery of genuine multivalued functions satisfying $f_1+ \cdots + f_N = 1$,
but does not yet incorporate the nonnegativity constraints $f_n \ge 0$.
To this end,
one considers more generally model sets of the type $\cK \cap \cC$,
i.e., intersections of a convex model set $\cK$ 
(e.g. an approximability set) with a convex cone $\cC$
(e.g. the set of nonnegative functions).
Thanks to Lemma \ref{Lem2},
as soon as $\cK$ and $\cC$ are amenable to computations,
the program \eqref{MainOptAbs}
turns into the manageable convex program presented below.

\bthm
\label{ThmIntersection}
For model sets of the form $\cK_n \cap \cC_n$,
where both $\cK_n$ and $\cC_n$ are convex and $\cC_n$ is a cone,
the affine optimal recovery map of Theorem~\ref{ThmAbs} is constructed by solving, for each $j \in [1:N]$,
the constrained convex optimization program \vspace{-2mm}
\begin{align*}
\minimize{\substack{ \bm{\mu}^+, \bm{\mu}^- \in (F^*)^L\\ \bm{\nu}^+, \bm{\nu}^- \in (F^*)^N \\ c \in \bR^{M \times N}, \, d,e \in \bR }} 
\; \,  e 
& & \mbox{s.to} & \quad
\pm \sum_{\ell=1}^L \mu^\pm_\ell(b_\ell)
+ \sum_{n=1}^N \tn{\nu^\pm_n}_{\cK_n} \le e \pm d\\
& & \mbox{and} & \qquad
\nu^\pm_n \mp 
\Big(
\delta_{j,n} \gamma_j - \sum_{m=1}^M c_{m,n} \la_m - \sum_{\ell=1}^L a_{\ell,n} \mu^\pm_\ell \Big)
\in \cC_n^{\rm dual}
\quad \mbox{for all }n\in [1:N].
\end{align*}
\ethm

\bpf
For each $\eta^\pm_n := \pm \delta_{j,n} \gamma_j \mp \sum_{m=1}^M c_{m,n} \la_m - \sum_{\ell=1}^L a_{\ell,n} \mu^\pm_\ell$
appearing in Theorem \ref{ThmAbs},
the support function of $\cK_n \cap \cC_n$ evaluated at $\eta^{\pm}_n$ is transformed, according to Lemma \ref{Lem2},
into
\begin{align*}
\tn{ \eta^\pm_n}_{\cK_n \cap \cC_n}
& = \sup_{f_n \in \cK_n \cap \cC_n} \eta^\pm_n (f_n)
= \min_{\nu^\pm_n \in F^*}
\Big\{ \sup_{f_n \in \cK_n} \nu^\pm_n(f_n): 
\, \nu^\pm_n - \eta^\pm_n \in \cC_n^{\rm dual} \Big\}\\
& = \min_{\nu^\pm_n \in F^*} \big\{ \tn{\nu^\pm_n}_{\cK_n}:
\, \nu^\pm_n - \eta^\pm_n \in \cC_n^{\rm dual} \big\}.
\end{align*}
Incorporating the $\nu^\pm_n$'s as optimization variables in \eqref{MainOptAbs}
while also making the change $\bm{\mu}^- \leftrightarrow -\bm{\mu}^-$ leads to the announced optimization program.
\epf

\section{Instantiation of the Main Results}
\label{SecInst}
\vspace{-3mm}

In this subsection, the abstract formalism adopted so far is 
specified in two situations of practical interest, 
leading to the two theorems highlighted in the introduction.
These situations are that of a (reproducing kernel) Hilbert space and of  a space of continuous functions.

\subsection{Hilbert spaces}
\label{SubSecHil}
\vspace{-2mm}

One assumes in this subsection that $F$ is a Hilbert space---hence the notation $H$ is used instead.
In this situation, any linear functional $\mu \in H^*$
is identified with its Riesz representer $u \in H$ via
$\mu(f) = \langle u, f \rangle$ for all $f \in H$.
In particular, the linear functionals $\la_m \in H^*$ and $\gamma_j \in H^*$ are identified with vectors $w_m \in H$ and $g_j \in H$.
Thus,
when $\la_m = \delta_{x^{(m)}}$ and $\gamma_j = \delta_{x^{(0)}}$ are point evaluations in a reproducing kernel Hilbert space with kernel~$K$,
one has $w_m = K(x^{(m)},\cdot)$ and $g_j = K(x^{(0)},\cdot)$.
Furthermore, the determination \eqref{SupFunAppSet} of the support function of the set $\cK_{V,T}$ defined in \eqref{GenApproxSets},
evaluated at a linear functional identified with $u \in H$,
easily reduces to
$$
\tn{u}_{\cK_{V,T}} =
\begin{cases}
\| T^{-*} u\|_H & \text{if } T^{-*} u \in V^\perp,\\
+\infty & \text{otherwise.} 
\end{cases}
$$
Taking this observation into account,
Theorem \ref{ThmAbs4AppSet} yields the corollary below,
of which Theorem \ref{ThmIntroHilb} is the special case $T_n = (1/\eps_n) \Id_H$, $L=1$, $A= [1;\ldots;1]$, and $b = 1$
(again assuming that $H$ is a reproducing kernel Hilbert space containing the constant function $1$).

\bcor
\label{CorResHil}
In a Hilbert space $H$,
if the model sets are $\cK_{V_n,T_n} = \{ f \in H: {\rm dist}_H(T_n f,V_n) \le 1 \}$,
then the affine optimal recovery map of Theorem~\ref{ThmAbs} is constructed by solving, for each $j\in[1:N]$,
the convex optimization program \vspace{-2mm}
\begin{align}
\nonumber
\minimize{\substack{ \vu^+, \vu^- \in H^L \\ c \in \bR^{M \times N}, \, d,e \in \bR }} 
 \;   e 
 & & \mbox{s.to} & \quad 
\pm \sum_{\ell = 1}^L \langle u^\pm_\ell, b_\ell \rangle
+ \sum_{n=1}^N 
\Big\|T_n^{-*} \Big( \delta_{j,n} g_j - \sum_{m=1}^M c_{m,n} w_m - \sum_{\ell=1}^L a_{\ell,n} u^\pm_\ell \Big) \Big\|_H
\le e \pm d \\
\label{ProgramHilbert}
& & \mbox{and} & \quad
\Big\langle T_n^{-*} \Big( \delta_{j,n} g_j - \sum_{m=1}^M c_{m,n} w_m - \sum_{\ell=1}^L a_{\ell,n} u^\pm_\ell \Big),
v \Big\rangle = 0 \;  \mbox{for all }v \in V_n, n \in [1:N].
\end{align}
\ecor

Should the norms be squared in the above,
the optimization program would have a closed-form solution.
Such a favorable outcome cannot be hoped for here,
since \eqref{ProgramHilbert} includes the geometric median problem as a special case 
and the latter is not known to possess closed-form solutions.
Regardless, there is no difficulty in solving \eqref{ProgramHilbert} computationally---it is a second-order cone program.
However, it is not immediately clear how to handle an infinite-dimensional Hilbert space $H$,
but fortunately a representer theorem holds,
at least in some cases.
Precisely,
it is possible to replace the minimization over $H$ by a minimization over  a finite-dimensional subspace of $H$ when the $T_n$'s are of the form $T_n = (1/\eps_n)T$ for a common operator $T$.
To justify this,
it is enough to remark that if $u_1^\pm,\ldots,u_L^\pm,c,d,e$ are feasible for  
\eqref{ProgramHilbert},
then so are $\wt{u}_1^\pm, \ldots, \wt{u}_L^\pm,c,d,e$,
where $\wt{u}_\ell^\pm := T^*P_{\wt{H}}T^{-*}u_\ell^\pm$
for $\ell \in [1:L]$ and where 
$P_{\wt{H}}$ denotes the orthogonal projector onto the finite-dimensional space
$$
\wt{H} := {\rm span}\{T b_1,\ldots, T b_L, T^{-*} g_j, T^{-*}w_1, \ldots, T^{-*}w_M \} + V_1 + \cdots + V_N .
$$
This remark follows from the fact that, since the $T b_\ell$'s belong to  $\wt{H}$,
$$
\langle \wt{u}^\pm_\ell, b_\ell \rangle
= \langle T^{-*} u^\pm_\ell, P_{\wt{H}} T b_\ell \rangle
= \langle T^{-*} u^\pm_\ell, T b_\ell \rangle
= \langle u^\pm_\ell, b_\ell \rangle,
$$
from the fact that, since $T^{-*} g_j$ and the $T^{-*}w_m$'s belong to $\wt{H}$,
\begin{align*}
\Big\| & T_n^{-*}  \Big( \delta_{j,n} g_j - \hspace{-1mm} \sum_{m=1}^M c_{m,n} w_m  - \hspace{-1mm} \sum_{\ell=1}^L a_{\ell,n} \wt{u}^\pm_\ell \Big) \Big\|_H = 
\eps_n \Big\| T^{-*} \Big( \delta_{j,n} g_j - \hspace{-1mm} \sum_{m=1}^M c_{m,n} w_m \Big) - \hspace{-1mm} \sum_{\ell=1}^L a_{\ell,n} P_{\wt{H}} T^{-*} u^\pm_\ell \Big\|_H\\
& \hspace{-1mm} = \eps_n \Big\| P_{\wt{H}} T^{-*} \big( \delta_{j,n} g_j - \hspace{-1mm} \sum_{m=1}^M c_{m,n} w_m - \hspace{-1mm} \sum_{\ell=1}^L a_{\ell,n} u^\pm_\ell \big) \Big\|_H
\le \eps_n \Big\| T^{-*} \big( \delta_{j,n} g_j - \hspace{-1mm} \sum_{m=1}^M c_{m,n} w_m - \hspace{-1mm} \sum_{\ell=1}^L a_{\ell,n} u^\pm_\ell \big) \Big\|_H \\
& \hspace{-1mm} =
\Big\|T_n^{-*} \Big( \delta_{j,n} g_j - \hspace{-1mm} \sum_{m=1}^M c_{m,n} w_m - \hspace{-1mm} \sum_{\ell=1}^L a_{\ell,n} u^\pm_\ell \Big) \Big\|_H,
\end{align*} 
and from the fact that, since all $v \in V_n$, $n \in [1:N]$, belong to $\wt{H}$,
\begin{align*}
\Big\langle & T_n^{-*} \Big( \delta_{j,n} g_j - \hspace{-1mm} \sum_{m=1}^M c_{m,n} w_m - \hspace{-1mm} \sum_{\ell=1}^L a_{\ell,n} \wt{u}^\pm_\ell \Big),
v \Big\rangle
 = \eps_n
 \Big\langle P_{\wt{H}} T^{-*} \Big( \delta_{j,n} g_j - \hspace{-1mm} \sum_{m=1}^M c_{m,n} w_m - \hspace{-1mm} \sum_{\ell=1}^L a_{\ell,n} u^\pm_\ell \Big),
v \Big\rangle\\
& = \eps_n \Big\langle T^{-*} \Big( \delta_{j,n} g_j - \hspace{-1mm} \sum_{m=1}^M c_{m,n} w_m - \hspace{-1mm} \sum_{\ell=1}^L a_{\ell,n} u^\pm_\ell \Big),
v \Big\rangle
=  \Big\langle T_n^{-*} \Big( \delta_{j,n} g_j - \hspace{-1mm} \sum_{m=1}^M c_{m,n} w_m - \hspace{-1mm} \sum_{\ell=1}^L a_{\ell,n} u^\pm_\ell \Big),
v \Big\rangle.
\end{align*}

\subsection{Spaces of continuous functions}
\label{SubSecCX}
\vspace{-2mm}

One assumes in this subsection that $F$ is the space $C(\cX)$ of continuous functions on a compact set~$\cX$.
In this situation,
the linear functionals $\mu \in C(\cX)^*$ are identified---keeping the same notation---with signed Borel measures $\mu \in B(\cX)$ via $\mu(f) = \int_\cX f \, d \mu$ for all $f \in C(\cX)$.
The attention is put on approximability sets intersected with the cone $\cC$ of nonnegative functions.
Here, the crucial ingredients are the facts that
the norm of $\mu \in C(\cX)^*$
is the total variation $\int_\cX d |\mu|$ of $\mu \in B(\cX)$
and that the dual cone $\cC^{\rm dual}$ of $\cC$ is the set $B_+(\cX)$ of nonnegative Borel measures.
Thus, Theorem~\ref{ThmIntersection} yields the following corollary,
of which Theorem \ref{ThmIntroCont} is the special case $L=1$, $A=[1;\ldots;1]$, $b=1$.

\bcor
\label{CorCNonneg}
In a space $C(\cX)$ of continuous functions in a compact set $\cX$,
if the model sets are $\cK_n = \{ f \in C(\cX): {\rm dist}_{C(\cX)}(f,V_n) \le \eps_n, f \ge 0 \}$,
then the affine optimal recovery map of Theorem~\ref{ThmAbs} is constructed by solving, for each $j \in [1:N]$,
the convex optimization program  \vspace{-2mm}
\begin{align*}
\minimize{\substack{ \bm{\mu}^+, \bm{\mu}^- \in B(\cX)^L\\ \bm{\nu}^+, \bm{\nu}^- \in B(\cX)^N \\ c \in \bR^{M \times N}, \, d,e \in \bR }} 
\; \,  e 
& & \mbox{s.to} & \;
\pm \sum_{\ell=1}^L \int_{\cX} b_\ell \, d\mu^\pm_\ell
+ \sum_{n=1}^N \eps_n \int_\cX d|\nu^\pm_n|
\le e \pm d\\
& & \mbox{and} & \;
\int_\cX v \, d{\nu^\pm_n} = 0 \quad
\mbox{for all } v \in V_n, \; n \in [1:N],\\
& & \mbox{and} & \; 
\nu^\pm_n \mp \Big( \delta_{j,n} \gamma_j - \sum_{m=1}^M c_{m,n} \la_m - \sum_{\ell=1}^L a_{\ell,n} \mu^\pm_\ell \Big)
\in B_+(\cX),
\quad n \in [1:N].
\end{align*}
\ecor

Despite the convexity of the above program,
its infinite dimensionality calls for a justification of its computational implementation.
In a nutshell,
it is a minimization over measures 
and as such it can be attacked via the Moment-SoS hierarchy
(see \cite{JBL} for a \mbox{survey}).
The guiding arguments are sketched below in the simpler univariate setting, e.g. considering each $f_n(x)$ as the concentration of the $n$th constituent as time $x$ evolves in $[0,\pi]$, say.
For convenience,
it is assumed that each space~$V_n \in C[0,\pi]$ is made of cosine polynomials of degree at most $K(n)$
and that each $b_\ell \in C[0,\pi]$ is also a cosine polynomial.

First,
one decomposes each optimization variable $\eta \in B(\cX)$ as the difference of its positive part $\eta^{\oplus}\in B_+(\cX)$ and its negative part $\eta^{\ominus} \in B_+(\cX)$,
so that $\eta = \eta^{\oplus}-\eta^{\ominus}$ and $|\eta| = \eta^{\oplus} + \eta^{\ominus}$.
Next, one thinks of the nonnegative Borel measure $\eta^{\tiny \mbox{\textcircled{$\pm$}}}$ equivalently in terms of its (infinite) sequence $y^{\tiny \mbox{\textcircled{$\pm$}}} \in \bR^\bN$ of trigonometric moments defined by $y^{\tiny \mbox{\textcircled{$\pm$}}}_k = \int_0^\pi \cos(k \theta) d \eta^{\tiny \mbox{\textcircled{$\pm$}}}(\theta)$,
$k \ge 0$, 
provided that the (infinite) symmetric Toeplitz matrix
${\rm Toep}_\infty(y^{\tiny \mbox{\textcircled{$\pm$}}})$
 built from it is positive semidefinite (see \cite[Section 3]{FouLas} for details).
Then, expressing each $b_\ell \in C[0,\pi]$ 
as $b_\ell(\theta) = \sum_{k} b_{\ell,k} \cos(k \theta)$---this is particularly neat for the dependence relation $f_1 + \cdots + f_N = 1$---one arrives at the equivalent infinite-dimensional semidefinite program \vspace{-2mm}
\begin{align}
\nonumber
\minimize{\substack{ \vu^{+,\tiny \mbox{\textcircled{$\pm$}}},
\vu^{-,\tiny \mbox{\textcircled{$\pm$}}} \in (\bR^\bN)^L 
\\ \vv^{+,\tiny \mbox{\textcircled{$\pm$}}},
\vv^{-,\tiny \mbox{\textcircled{$\pm$}}} \in (\bR^\bN)^N \\ c \in \bR^{M \times N}, \, d,e \in \bR }} 
\; \,  e 
& & \mbox{s.to} & \;
\pm \sum_{\ell=1}^L \sum_k  b_{k,\ell}  (u^{\pm,\oplus}_{\ell,k} - 
u^{\pm,\ominus}_{\ell,k}) 
+ \sum_{n=1}^N \eps_n (v^{\pm,\oplus}_{n,0} + v^{\pm,\ominus}_{n,0})
\le e \pm d\\
\nonumber
& & \mbox{and} & \;
v^{\pm,\oplus}_{n,k} - v^{\pm,\ominus}_{n,k} = 0 \quad
\mbox{for all } k \in [0:K(n)], \; n \in [1:N],\\
\nonumber
& & \mbox{and} & \; 
{\rm Toep}_\infty\Big(
v^{\pm,\oplus}_n - v^{\pm,\ominus}_n \mp \delta_{j,n} y_j
\pm \sum_{m=1}^M c_{m,n} z_m
\pm \sum_{\ell=1}^L a_{\ell,n} (u^{\pm,\oplus}_\ell - u^{\pm,\ominus}_\ell)
\Big) \succeq 0,\\
\label{SDPInf}
& & \mbox{and} & \; 
{\rm Toep}_\infty\big( u^{\pm,\tiny \mbox{\textcircled{$\pm$}}}_\ell \big) \succeq 0,
\quad \ell \in [1:L],
\qquad 
{\rm Toep}_\infty\big( v^{\pm,\tiny \mbox{\textcircled{$\pm$}}}_n \big) \succeq 0,
\quad n \in [1:N].
\end{align}
Above, the fixed  $y_j,z_1,\ldots,z_M \in \bR^\bN$ represented the (infinite) sequences of trigonometric moments stemmed from $\gamma_j, \la_1,\ldots, \la_M$.
In the case $\gamma_j = \delta_{x^{(0)}}$ and $\la_m = \delta_{x^{(m)}}$,
these moments are given by $y_{j,k} = \cos(k x^{(0)})$
and $z_{m,k} = \cos(k x^{(m)})$,
$k \ge 0$.
Finally, in order to deal with manageable semidefinite programs,
one truncates the optimization variables $u^{\pm,\tiny \mbox{\textcircled{$\pm$}}}_\ell, v^{\pm,\tiny \mbox{\textcircled{$\pm$}}}_n$ at a level $r$,
so that they belong to $\bR^{r}$ rather than $\bR^\bN$,
and one concurrently replaces the infinite Toeplitz matrices 
by their $r \times r$ upper-left sections.
It can be shown that the minimizers of these truncated programs converge to minimizers of the infinite-dimensional program as $r \to \infty$.
This is the essence of the semidefinite hierarchy.

But even without being ensured of the convergence,
the above procedure supplies a lower bound ${\rm lb}_r$ for the optimal value of the original program
at each level $r$,
since truncating an infinite sequence feasible for \eqref{SDPInf}
yields a finite sequence feasible for the truncated program.
In other words,
one can compute a value ${\rm lb}_r \le \inf \{ {\rm gwce}_\infty(\Delta), \, \Delta: \bR^{M \times N} \to \bR^N \}$. 
The associated minimizers $c^{[r,j]} \in \bR^{M \times N}$ and $d^{[r,j]} \in \bR$ allow one to construct an affine recovery map $\Delta^{[r]} = [\Delta^{[r]}_1;\ldots;\Delta^{[r]}_N]$ via 
$\Delta^{[r]}_j: y \in \bR^{M \times N} \mapsto \sum_{m=1}^M \sum_{n=1}^N c_{m,n}^{[r,j]} y_{m,n} + d^{[r,j]}$.
Its worst-case error is the maximal optimal value of $N$ convex optimization programs involving measures (stemming from the proof of Theorem \ref{ThmAbs} and left for the reader to spell out). 
An upper bound  ${\rm ub}_{r,s}$ on ${\rm gwce}_\infty(\Delta^{[r]})$ can then be obtained by imposing these measures to be atomic measures on an $s$-point grid---hence turning the convex programs into a linear programs.
In other words, one can compute a value ${\rm ub}_{r,s} \ge {\rm gwce}(\Delta^{[r]})$.
All in all, one can now certify that 
$$
{\rm gwce}(\Delta^{[r]}) \le \f{{\rm ub}_{r,s}}{{\rm lb}_{r}} \,
\inf \{ {\rm gwce}_\infty(\Delta), \, \Delta: \bR^{M \times N} \to \bR^N \},
$$
i.e., that $\Delta^{[r]}$ is a near-optimal recovery map with factor ${\rm ub}_{r,s}/{\rm lb}_{r}$
that can be estimated {\em a posteriori}
(and that should approach one as $r$ and $s$ grow).

\section{Parting Thoughts}
\label{SecConc}
\vspace{-3mm}

As far as the author is aware,
the practical reconstruction of multivalued functions was not treated in the Optimal Recovery literature before.
This article starts to fill the gap,
but there are several open questions still to be investigated,
some of them discussed below.\vspace{-2mm}

{\bf Other measures of performance.}
The main results of this article dealt with worst-case errors using the $\ell_\infty^N$-norm to quantify the prediction of $\vf(x^{(0)}) \in \bR^N$.
What about $\ell_p^N$ for other $p \in [1,\infty)$, say for $p=1$?
By virtue of the inequalities $\|z\|_{\ell_\infty^N} \le \|z\|_{\ell_1^N} \le N \|z\|_{\ell_\infty^N}$ for $z \in \bR^N$,
one easily sees that the affine optimal recovery map relative to $\ell_\infty^N$
is near-optimal relative to $\ell_1^N$ with factor~$N$.
But this can be improved:
according to Proposition (i) of \ref{PropApp} in the appendix,
and in view of the estimation ${\rm proj}(\ell_1^N) \lesssim \sqrt{2/\pi} \sqrt{N}$ for the projection constant of $\ell_1^N$ originally due to \cite{Gru},
the near-optimality factor for affine maps is of order at most $\sqrt{N}$.
Can this be further reduced, say to a constant?
Importantly, can one construct a genuinely optimal recovery map, affine or not? \vspace{-2mm}

{\bf Full Recovery.}
The article mostly tackled the prediction of multivalued functions at a fixed~$x$,
i.e.,
the recovery of $\Gamma = \delta_x$.
The ultimate goal would be the recovery of $\Gamma = \Id_F$.
It can be seen that assembling affine optimal maps for the recovery of each $\delta_x$, $x \in \cX$,
provides an affine optimal map for the full recovery of $\Gamma = \Id_F$ if one works within $F=L_\infty(\cX,\ell_\infty^N)$.
It is conceivable that this assembly remains optimal within
$F=C(\cX,\ell_\infty^N)$---this was established, not without efforts,
for $N=1$ and approximability sets in \cite{DFPW}.
Such an assembly is not a practical construction, though,
and more efforts will be needed to achieve one---in the case $N=1$,
this was done for specific approximability sets in \cite{FullChe}.
As for full recovery within a Hilbert space $H$,
the problem seems even more intricate,
because it includes (with the choice $A \vf = [f_1-f_2;\cdots;f_N - f_{N-1}]$ and $\vb = 0$)
the full recovery problem relative to a model set equal to the intersection of $N$ hyperellipsoids,
which has only been solved for $N=2$ in \cite{2Hyp}. \vspace{-2mm}

{\bf Observation Errors.}
The article only studied exact observations of the type $y_{m,n} = f_n(x^{(m)})$,
but in realistic situations,
these observations are inaccurate,
i.e., of the type $y_{m,n} = f_n(x^{(m)}) + e_{n,m}$.
The error vectors $e_1,\ldots,e_N \in \bR^M$ can be modeled stochastically as random vectors.
In this case, considering the prediction problem,
one can optimistically hope for near-optimality of affine maps,
as obtained in \cite{Don} for Gaussian noise and in \cite{FouPao} even for log-concave noise with an altered measure of performance.
The errors vectors $e_1,\ldots,e_N \in \bR^M$ can also be modeled deterministically as belonging to uncertainty sets $\cE_1,\ldots,\cE_N$.
In this case, there is a classical reduction to the accurate scenario by introducing the compound objects $[f_n;e_n] \in F \times \bR^M$,
about which one has the {\em a priori} information $[f_n;e_n] \in \cK_n \times \cE_n$ and $A \vf = \vb$,
as well as the {\em a posteriori} information $\wt{\la}([\vf;\ve]) := \La(\vf) + \ve$,
and one aims at recovering the quantity of interest $\wt{\Gamma}([\vf;\ve]) := \Gamma (\vf)$.
The abstract theory developed here provides constructions of optimal recovery maps,
but may require additional work to turn them into practical constructions,
for instance when the model and uncertainty sets are not independent---e.g., when the functions $f_n$ represent concentrations, the observations $y_{m,n}$ are likely to be measured as concentrations, too,
hence assumptions $0 \le f_n(x) \le 1$ and $|e_{n,m}| \le \eta$, say,
need to be augmented with the $0 \le f_n(x^{(m)}) + e_{n,m} \le 1$.

\section*{Appendix}
\vspace{-3mm}

This section goes back one step and reconsiders the general Optimal Recovery framework described in Subsection \ref{SSecOverviewOR}.
Proposition \ref{PropApp} below amalgamates,
through one fresh take,
some old and new results about near-optimality of affine recovery maps relative to convex model sets.
Since it involves the notion of projection constants,
one recalls that the relative projection constant of a subspace~$V$ of a norm space $W$ is defined by
$$
{\rm proj}(V,W) := \inf \big\{
\|P\|: \; P \mbox{ is a projection from $W$ onto $V$}
\big\},
$$
where $P$ being a projection from $W$ onto $V$ means that $P$ is a linear map from $W$ into $V$ such that $Pv= v$ for all $v \in V$.
As for the absolute projection constant of a normed space $U$,
it is defined by
$$
{\rm proj}(U) := \sup \{ {\rm proj}(i(U),W), \; i \mbox{ is an isometric embedding from $U$ into $W$} \}.
$$ 
The absolute projection constant equals the extension constant (see e.g. \cite[Theorem~5]{Woj}), i.e.,
\begin{align*}
{\rm proj}(U) = \inf \big\{
c: & \mbox{ for any }V \inc W \mbox{ and any linear map } T:V \to U,\\
& \mbox{ there exists a linear map } \wt{T}: W \to U \mbox{ such that }
T_{|V} = T \mbox{ and } \|\wt{T}\| \le c \|T\| \big\}.
\end{align*}
To retrieve the results stated in the text, one needs  Kadec--Snobar estimate ${\rm proj}(U) \le \sqrt{\dim(U)}$ 
(see \cite[Theorem 10]{Woj}) and
its corollary ${\rm proj}(U) \le 1 + \sqrt{{\rm codim}(U)}$
(see \cite[Corollary~11]{Woj}).

\bprop
\label{PropApp}
Given a vector space $F$ and a model set $\cK \inc F$,
given a linear observation map $\La: F \to \bR^M$,
and given a linear quantity of interest $\Gamma: F \to Z$,\vspace{-5mm}
\begin{enumerate}
\item[(i)]
If $\cK$ is convex, \vspace{-2mm}
\begin{enumerate}
\item[(i-a)] then there exists an affine recovery map which is near optimal with factor $2C$, where $C:= \min\{ {\rm proj}(Z), 1 + \sqrt{M} \}$;\vspace{-1mm}
\item[(i-b)] and if $\Gamma: F \to \bR$ is furthermore a linear functional,
then there exists an affine recovery map which is genuinely optimal. \vspace{-2mm}
\end{enumerate}
\item[(ii)]
If $\cK$ is convex and contains the origin, then affine can be replaced by linear. \vspace{-2mm}
\item[(iii)]
If $\cK$ is convex and symmetric about the origin,
then the factor $2C$ can be replaced by $C$. \vspace{-2mm}
\end{enumerate} 
\eprop

\bpf
Central to the forthcoming arguments is the lower bound on the worst-case error of a recovery map $\Delta:\bR^M \to Z$ by half of the so-called diameter of information, i.e.,
$$
{\rm gwce}(\Delta) := \sup_{f \in \cK} \|\Gamma(f) - \Delta(\La f) \|_Z 
\ge {\rm lb} := \sup_{\substack{ h \in \ker(\La) \\ h \in (\cK - \cK)/2 }} \|\Gamma(h)\|_Z.
$$
This can be easily justified in exactly the same way as in the beginning of the proof of Proposition~\ref{PropIndep}.
Since $\cK$ is convex,
the set $\cK_{\rm sym}: = (\cK - \cK)/2$ is convex and symmetric about the origin,
so its Minkowski functional is a seminorm,
denoted by $| \cdot |$ here. 
The defining expression of the lower bound yields
\be
\label{beforeExt}
\|\Gamma(h)\|_Z \le {\rm lb} \, |h|
\qquad
\mbox{for all } h \in \ker(\La).
\ee
It is implicitly assumed that the lower bound is a finite quantity for any $\Gamma$,
in particular for $\Gamma = \Id_F$,
which implies that $\ker(\La) \cap G' = \{0\}$
where $G' := \{ g' \in F: |g'| = 0\}$ is a linear subspace of $F$.
Now consider another linear subspace $G''$ containing $\ker(\La)$  such that $G' \oplus G'' = F$
and notice that $| \cdot |$ induces a norm on $G''$.
On the one hand,
by the definition of the extension constant of $Z$,
there exists a linear map $\wh{\Gamma}: G'' \to Z$ such that $\wh{\Gamma}_{| \ker(\La)} = \Gamma_{| \ker(\La)}$ and
$$
\|\wh{\Gamma}(g'')\|_Z \le {\rm proj}(Z) \, {\rm lb} \, |g''|
\qquad
\mbox{for all } g'' \in G''.
$$
On the other hand,
by selecting a projection $P$ from $G''$ onto $\ker(\La)$ with operator norm at most $ 1+\sqrt{M}$
and by defining a linear map $\wh{\Gamma}: G'' \to Z$ via
$\wh{\Gamma}(g'') = \Gamma(P(g''))$ for $g'' \in G''$,
one obtains
$$
\|\wh{\Gamma}(g'')\|_Z 
\le {\rm lb} \, |P(g'')|
\le (1+\sqrt{M}) \, {\rm lb} \, |g''| 
\qquad
\mbox{for all } g'' \in G''.
$$
Therefore, one always has $\|\wh{\Gamma}(g'')\|_Z  \le C \, {\rm lb} \, |g''|$ for all $g'' \in G''$ for some $\wh{\Gamma}: G'' \to Z$ satisfying $\wh{\Gamma}_{| \ker(\La)} = \Gamma_{| \ker(\La)}$.
Next, one introduces the linear map $\wt{\Gamma}: F \to Z$ defined by
$\wt{\Gamma}(g') = 0$ for $g' \in G'$
and $\wt{\Gamma}(g'') = \wh{\Gamma}(g'')$ for $g'' \in G''$.
Noticing that $\wt{\Gamma}(g) = \wh{\Gamma}(g'')$
and $|g| = |g''|$ for any $g \in F$ written as $g = g' + g''$ with $g' \in G'$ and $g'' \in G''$,
one arrives at  
$$
\| \wt{\Gamma} (g) \|_Z
\le C \, {\rm lb} \, |g|
\qquad
\mbox{for all } g \in F.
$$
It also holds that $(\Gamma - \wt{\Gamma})_{| \ker(\La)} = 0$, which implies the existence of $c_1,\ldots,c_M \in Z$ such that 
$\Gamma - \wt{\Gamma} = \sum_{m=1}^M c_m \la_m$.
Thus, the previous estimate reads
\be
\label{KeyBound}
\Big\| \Big(  \Gamma - \sum_{m=1}^M c_m \la_m \Big)(g) \Big\|
\le C \, {\rm lb} \, |g|
\qquad
\mbox{for all } g \in F.
\ee
Now, given $f \in \cK$ and
picking an arbitrary $\ol{f} \in \cK$,
applying the above to  $g := (f-\ol{f})/2 \in \cK_{\rm sym}$, 
which satisfies $| g | \le 1$, yields
$$
\Big\| \Big(  \Gamma - \sum_{m=1}^M c_m \la_m \Big) f
- \Big(  \Gamma - \sum_{m=1}^M c_m \la_m \Big) \ol{f} \Big\|
\le 2 C \, {\rm lb},
$$
i.e., $\|\Gamma(f) - \Delta^{\rm aff}(\La f)\|_Z \le 2C \, {\rm lb}$,
where the affine recovery map $\Delta^{\rm aff}$ is defined for $y \in \bR^M$ by
$\Delta^{\rm aff}(y) = \sum_{m=1}^M c_m y_m + d$
with $d:=  \Gamma(\ol{f}) - \sum_{m=1}^M c_m \la_m (\ol{f})$.
Taking the supremum over $f \in \cK$
reveals that ${\rm gwce}(\Delta^{\rm aff}) \le 2C \, {\rm lb} \le 2C \, \inf\{ {\rm gwce}(\Delta), \Delta: \bR^M \to Z \}$,
i.e., that the recovery map $\Delta^{\rm aff}$ is near optimal with factor $2C$.
This establishes (i-a).
For (ii), note that,
if $0 \in \cK$, then one can pick $\ol{f} = 0$ to turn $\Delta^{\rm aff}$ into a linear map. 
For (iii), 
note that, if $\cK$ is symmetric about the origin,
then $\cK_{\rm sym}$ coincides with $\cK$ and one can apply \eqref{KeyBound} directly to $ g = f \in \cK$, so that $|g| \le 1$,
and derive near optimality with factor $C$ instead of $2C$.

It remains to justify (i-b) when $\Gamma$ is a linear functional---this is the result of \cite{Suk},
proved in a different manner here.
One comes back to \eqref{beforeExt} written for $Z = \bR$,
i.e.,
$$
\Gamma(h) \le {\rm lb} \, |h|
\qquad
\mbox{for all } h \in \ker(\La).
$$
The dominated extension theorem guarantees the existence of an extension $\wt{\Gamma}$ to $\Gamma_{|\ker(\La)}$ to the whole $F$ such that $$
\wt{\Gamma}(g) 
\le {\rm lb} \, |g|
\qquad
\mbox{for all } g \in F.
$$
Therefore, for all $f^+,f^- \in \cK$,
one has  $(\wt{\Gamma}(f^+) - \wt{\Gamma}(f^-))/2 \le {\rm lb}$,
i.e., $(S-I)/2 \le {\rm lb}$,
where $S:= \sup \{ \wt{\Gamma}(f), f \in \cK \}$ 
and $I := \inf \{ \wt{\Gamma}(f), f \in \cK \}$.
From $\wt{\Gamma}(f) \in [I,S]$ for all $f \in \cK$,
one obtains $ \wt{\Gamma}(f) - d  \in [-(S-I)/2, (S-I)/2] \inc [- {\rm lb}, {\rm lb}]$ where $d:= (I+S)/2$.
Since the extension $\wt{\Gamma}$ can be written as $\wt{\Gamma} = \Gamma - \sum_{m=1}^M c_m \la_m$ for some $c_1,\ldots,c_M \in \bR$,
one derives
$$ 
\Big|  \Gamma(f) - \Big( \sum_{m=1}^M c_m \la_m(f) + d \Big) \Big| 
\le {\rm lb}. 
$$
Taking the supremum over $f \in \cK$ yields ${\rm gwce}(\Delta^{\rm aff}) \le {\rm lb} \le \inf\{ {\rm gwce}(\Delta), \Delta: \bR^M \to \bR \}$,
which shows that the affine recovery map $\Delta^{\rm aff}: y \in \bR^M \mapsto \sum_{m=1}^M c_m y_m + d$ is genuinely~optimal.
\epf

\end{document}